\documentclass{ifacconf}
\makeatletter
\let\old@ssect\@ssect 
\makeatother

\usepackage{graphics,graphicx,epsf,epstopdf}      
\usepackage{natbib}        
\usepackage{amsmath,amsxtra,amsfonts,amscd,amssymb,bm}
\usepackage{floatrow} 

\usepackage[normalem]{ulem} 
\usepackage{hyperref}	
\hypersetup{linktocpage,colorlinks,linkcolor=blue,anchorcolor=blue,citecolor=blue,urlcolor=blue}
\makeatletter
\def\@ssect#1#2#3#4#5#6{%
	\NR@gettitle{#6}
	\old@ssect{#1}{#2}{#3}{#4}{#5}{#6}
}
\makeatother
\edef\endfrontmatter{%
	\unexpanded\expandafter{\endfrontmatter}
	\noexpand\endNoHyper 
}

\usepackage{tikz,array}
\usetikzlibrary{shapes.geometric}
\usetikzlibrary{shapes.arrows}
\usetikzlibrary{calc}
\usetikzlibrary{intersections}

\newcommand{\Rnn}{\mathbb{R}^{n\times n}} 
\newcommand{\Rnp}{\mathbb{R}^{n\times p}} 
\newcommand{\Rpp}{\mathbb{R}^{p\times p}} 

\newcommand{\stiefel}{\mathrm{St}(p,n)}
\newcommand{\stM}{\mathrm{St}_{M}(p,n)}
\newcommand{\stX}{\mathrm{St}_{X^{\top}X}(p,n)}
\newcommand{\skewset}{{\cal S}_{\mathrm{skew}}}
\newcommand{\symset}{{\cal S}_{\mathrm{sym}}}
\newcommand{\TX}{{\mathrm{T}_{X}\stX}}
\newcommand{\TY}{{\mathrm{T}_{Y}\stM}}
\newcommand{\NX}{{\mathrm{N}_{X}\stX}}
\newcommand{\NY}{{\mathrm{N}_{Y}\stM}}
\newcommand{\proj}{\mathrm{P}_X}
\newcommand{\projn}{\proj^{\perp}}

\newcommand{\flow}{\varphi_t(X_0)}
\newcommand{\Zstar}{\mathcal{C}} 

\newcommand{\ff}{_{\mathrm{F}}}  
\newcommand{\fs}{^2_{\mathrm{F}}}
\newcommand{\inv}{^{-1}} 
\newcommand{\st}{\mathrm{s.\,t.}\,\,} 
\newcommand{\zz}{^{\top}} 
\newcommand{\dt}{\frac{\mathrm{d}}{\mathrm{d}t}}

\newcommand{\dkh}[1]{\left(#1\right)}

\newcommand{\hkh}[1]{\left\{#1\right\}}

\newcommand{\norm}[1]{\left\|#1\right\|}

\DeclareMathOperator*{\dimension}{dim}
\DeclareMathOperator*{\sym}{sym}
\DeclareMathOperator*{\skewsym}{skew}
\DeclareMathOperator*{\tr}{tr}

\usepackage[textwidth=\marginparwidth, textsize=tiny]{todonotes}

\newcommand{\pierre}[1]{\todo[color=blue!10, inline]{\textbf{Pierre:} #1}}

\newcommand{\bg}[1]{\todo[color=cyan!20, inline]{\textbf{Bin:} #1}}
\newcommand{\pacomm}[1]{\todo[color=magenta!10, inline]{\textbf{PA:} #1}}

\newcommand{\revise}[1]{#1}

\begin{document}
	\begin{frontmatter}
		
		\title{Optimization flows landing on the Stiefel manifold\thanksref{footnoteinfo}} 
		
		\thanks[footnoteinfo]{This work was supported by the Fonds de la Recherche Scientifique -- FNRS and the Fonds Wetenschappelijk Onderzoek -- Vlaanderen under EOS Project no. 30468160. 
			Bin Gao was supported by the Deutsche Forschungsgemeinschaft (DFG, German Research Foundation) via the collaborative research centre 1450--431460824, InSight, University of M\"unster, and via Germany's Excellence Strategy EXC 2044--390685587, Mathematics M\"unster: Dynamics--Geometry--Structure. Simon Vary is a beneficiary of the FSR Incoming Post-doctoral Fellowship.}
		
		\author[First]{Bin Gao}
		\author[Second]{Simon Vary}
		\author[Third]{Pierre Ablin} 
		\author[Second]{P.-A. Absil} 
		
		\address[First]{Institute for Applied Mathematics, University of M\"unster, 48149 Münster, Germany (e-mail: gaobin@lsec.cc.ac.cn)}
		\address[Second]{ICTEAM Institute, UCLouvain, 1348 Louvain-la-Neuve, Belgium (e-mail: \{simon.vary, pa.absil\}@uclouvain.be)}
		\address[Third]{CNRS, Université Paris-Dauphine, PSL University, France (e-mail: pierreablin@gmail.com).}
		
		\begin{abstract}                
			We study a continuous-time system that solves optimization problems over the set of orthonormal matrices, which is also known as the Stiefel manifold. The resulting optimization flow follows a path that is not always on the manifold but asymptotically lands on the manifold. We introduce a generalized Stiefel manifold to which we extend the canonical metric of the Stiefel manifold. We show that the vector field of the proposed flow can be interpreted as the sum of a Riemannian gradient on a generalized Stiefel manifold and a normal vector. Moreover, we prove that the proposed flow globally converges to the set of critical points, and any local minimum and isolated critical point is asymptotically stable.
		\end{abstract}
		
		
	\end{frontmatter}
	
	\section{Introduction}
	Consider the optimization problem
	\begin{equation}\label{prob:stiefel}
		\begin{array}{cl}
			\min\limits_{X\in\Rnp}&f(X)\\
			\st &  X\zz X = I_{p},
		\end{array}
	\end{equation}
	where $p\,{\le}\, n$, $I_p$ denotes the $p\times p$ identity matrix, the objective function $f:\Rnp\rightarrow\mathbb{R}$ is continuously differentiable, and the orthogonality constraints define the \emph{Stiefel manifold}, that is
	\begin{align*}
		\stiefel:=\hkh{X\in\Rnp: X\zz X = I_{p}}.
	\end{align*}
	Optimization over \revise{orthonormal} matrices as posed in~\eqref{prob:stiefel} appears in many practical applications, such as \revise{the} orthogonal procrustes problem~\citep{Elden1999Procrustes}, blind source \revise{separation}~\citep{Joho2002Joint}, \revise{the} linear eigenvalue problem~\citep{golub2013matrix},  principal component analysis~\citep{Grubisic2007Efficient} and its sparse variant~\citep{Chen2020Proximal}, electronic structure calculations~\citep{gao2020orthogonalization}; see \cite{edelman1998geometry} and \cite{Wen2013feasible} for a more complete list of applications. It is an instance of minimization over Riemannian manifolds for which many standard Euclidean algorithms have been extended \citep{AMS08, Hu2020Brief, boumal2022intromanifolds}. These methods are feasible, in that they follow a sequence of iterates that preserve the manifold constraint.
	
	
	Continuous-time systems have been used 
	{for solving matrix and optimization problems; see, e.g., \cite{Brockett1991dynamical}, \cite{Chu1994list}, \cite{MHM96}, \cite{Absil2006continuous}, and \cite{bournez2021survey}.}
	Recently, \cite{ablin2021fast} proposed a continuous-time flow called the {landing flow}, along with its discretization called the {landing algorithm}, that solves problem~\eqref{prob:stiefel} in the special case of the orthogonal manifold (that is, when $n\,{=}\,p$). The main advantage of the landing algorithm is that the individual iterates do not need to satisfy the manifold constraint, and therefore it alleviates the need to 
	{compute retractions that, depending on the objective function $f$, may be the computational bottleneck in optimization on the Stiefel manifold.}
	It is similar to an approach proposed in~\cite{gao2018parallelizable} for solving~\eqref{prob:stiefel}, which utilizes an augmented Lagrangian update that allows for a parallel implementation.
	
\begin{figure}[htbp]
		\centering
		\includegraphics[width=0.65\textwidth]{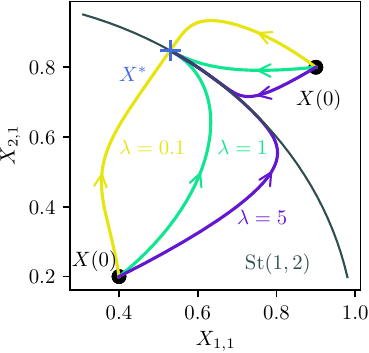}
		\caption{\revise{Landing flows on $\mathrm{St}(1,2)$ to minimize a linear function with different parameters~$\lambda$ and two {initial~points}.}\label{fig:diagram_landing}}
	\end{figure}

	In this paper, we extend the landing flow to solve the optimization problem~\eqref{prob:stiefel} over the Stiefel manifold; see Fig.~\ref{fig:diagram_landing} for an illustration. By considering a generalization of the Stiefel manifold and constructing a specific Riemannian metric, we give the landing flow a geometric interpretation \revise{involving a Riemannian gradient and a normal vector}. In addition, we prove the global convergence to the set of critical points  for the landing flow {and we study the stability of the equilibria}.
	
	
	This paper is organized as follows. 
	{After introducing the notation in~\S\ref{sec:notation}, we propose in~\S\ref{sec:landing} the landing flow on the}
	the Stiefel manifold. In~\S\ref{sec:interpretation}, a geometric interpretation of the landing flow is given. The convergence of the landing flow and  the stability of the equilibria are analyzed in~\S\ref{sec:convergence}. Finally, the conclusion is drawn in~\S\ref{sec:conclusion}.

	
	\section{Notation}\label{sec:notation}
	$\Rnp_*$ denotes the set of $n\times p$ matrices of full column rank. 
	{Given $X\,{\in}\,\Rnp_*$, we let $X_\perp$ denote an $n\times(n-p)$ matrix such that $X^\top X_\perp \,{=}\, 0$ and $X_\perp^\top X_\perp \,{=}\, I_{n-p}$.}
	The {Frobenius} inner product of two matrices $X, Y\,{\in}\, \Rnp$ is denoted by $\left\langle X,Y\right\rangle\,{=}\,\tr(X^{\top} Y)$, where $\tr(\cdot)$ denotes the matrix trace. The Frobenius norm of $X$ is denoted by $\|X\|\ff\,{:=}\,\sqrt{\left\langle X,X\right\rangle}$.
	$$\sym(A)=\frac12(A+A^{\top})\quad\mbox{and}\quad\skewsym(A)=\frac12(A-A^{\top})$$
	are the symmetric part and the skew-symmetric part of a \revise{square} matrix $A$, respectively.
	Moreover, $\symset^n$ and $\skewset^n$ denote the sets of all symmetric and skew-symmetric $n\times n$ matrices, respectively. 
	\revise{The Fr\'{e}chet derivative of {a~map} $F$ is denoted by $\mathrm{D}F$.} 
	
	\section{Landing flow on the Stiefel manifold} \label{sec:landing}
	{As an extension of the $p=n$ case addressed in~\cite{ablin2021fast}, we define the \emph{landing flow} on the Stiefel manifold as the solution}
	of a continuous-time system called the \emph{landing system}
	\begin{equation}\label{eq:landing-flow}
		\dot{X}(t) = - \Lambda\left(X\left(t\right)\right),
	\end{equation}
	where $\Lambda: \Rnp \rightarrow \Rnp$ is the \emph{landing field}\footnote{In order to give a geometric interpretation to $\Lambda(X)$ later, $\psi(X)$ differs from~\citet{ablin2021fast} by a factor \revise{of}~2.} defined as
	\begin{equation}\label{eq:landing-field}
		\Lambda(X) := \psi(X)X+\lambda \nabla\mathcal{N}(X).
	\end{equation}
	\revise{The} first component $\psi(X)X$ \revise{is} referred to as the \emph{relative gradient}, \revise{with}
	\begin{equation*}
		\psi(X):=2\skewsym\left(\nabla f(X)X^{\top}\right),
	\end{equation*}
	where $\nabla f(X)$ is the Euclidean gradient of the objective $f$. \revise{In the} second component of the landing field~\eqref{eq:landing-field},
	\begin{align*}
		\mathcal{N}(X):= \frac14 \left\|X^{\top} X -I_p\right\|\fs,
	\end{align*}
	whose set of minimizers is the Stiefel manifold, \revise{and $\lambda >0$ is a regularization parameter.}
	
	A crucial feature of the landing field is that we have $\nabla\mathcal{N}(X)=X(X\zz X-I_p)$ making the two components of the landing field in~\eqref{eq:landing-field} orthogonal
	\revise{with respect to the Frobenius inner product.}
	
	\cite{gao2018parallelizable} proposed a similar approach called PLAM for solving~\eqref{prob:stiefel} based on the field
	\begin{equation}\label{eq:PLAM-field}
		\tilde{\Lambda}(X) := \nabla f(X)-X\sym(\nabla f(X)^{\top} X ) +\lambda \nabla\mathcal{N}(X).
	\end{equation}
	{However, in~\eqref{eq:PLAM-field}, the component based on $\nabla f(X)$ is not orthogonal to $\nabla\mathcal{N}(X)$, contrary to the landing field~\eqref{eq:landing-field}.}
	
	Both fields, defined in~\eqref{eq:landing-field} and~\eqref{eq:PLAM-field}, have a straightforward interpretation as a Riemannian gradient for matrices $X$ belonging to the Stiefel manifold. To see this, notice that for $X\,{\in}\,\stiefel$, the distance term $\nabla\mathcal{N}(X)$ vanishes and, as a result, $\Lambda(X)$ and $\tilde{\Lambda}(X)$ belong to the tangent space of the Stiefel manifold
	$$ \mathrm{T}_X \stiefel = \{\xi\in\Rnp: \xi^{\top}X + X^{\top}\xi =0\}.$$
	Moreover, they are both Riemannian gradients of the objective $f$ on the Stiefel manifold
	but with respect to different Riemannian metrics. The landing field $\Lambda(X)$ corresponds to the Riemannian gradient with respect to the \emph{canonical metric} 
	\begin{equation}\label{eq:metric-canonical}
		g^{\mathrm{c}}_X(\xi,\zeta) := \langle \xi,(I_n - \frac12 XX^{\top})\zeta\rangle \quad \mbox{for~all}~\xi,\zeta\in\Rnp,
	\end{equation}
	while the PLAM field $\tilde{\Lambda}(X)$ corresponds to the Riemannian gradient with respect to the Euclidean metric $\langle \xi,\zeta\rangle$; see~\cite{edelman1998geometry} for the geometry of the Stiefel manifold and a discussion on these two different metrics.
	
	However, for a general matrix $X\,{\in}\,\Rnp_*$, which does not need to belong to the Stiefel manifold, the Riemannian interpretation of the landing field in~\eqref{eq:landing-field}~is not immediate. In the following section, we consider a certain generalization of the Stiefel manifold that allows us to derive a Riemannian gradient interpretation of the landing field for a general $X\,{\in}\,\Rnp_*$.
	
	\section{Interpretation of the landing flow} \label{sec:interpretation}
	We give a geometric interpretation of $\Lambda(X)\,{=}\,\psi(X)X+\lambda \nabla\mathcal{N}(X)$ for all $X\in\Rnp_*$ such that (i)~$\psi(X)X$ is the Riemannian gradient of $f$ on a Riemannian submanifold of $\Rnp_*$ with a specifically constructed metric and (ii)~$\nabla\mathcal{N}(X)$ belongs to its normal space. 
	
	First, we generalize the Stiefel manifold $\stiefel$ to the following set
	\begin{equation*} \label{eq:def_stM}
		\stM := \{Y \in \mathbb{R}^{n\times p}: Y^\top Y = M\},
	\end{equation*}
	where $M\,{\in}\,\Rpp$ is a given symmetric positive-definite matrix. In particular, $\stM$ reduces to the Stiefel manifold $\stiefel$ when $M=I_p$.
	
	Consider the linear map 
	\begin{equation}\label{eq:map-XtoY}
		\revise{\Phi_M}: \Rnp \to \Rnp: X \mapsto Y = X M^{\frac12}.
	\end{equation}
	\revise{It} is a diffeomorphism of $\Rnp_*$ onto itself since $M$ is symmetric positive definite, and it maps $\stiefel$ to $\stM$. As a consequence, the following proposition shows that $\stM$ is a submanifold of $\Rnp_*$.
	
	\begin{prop}
		$\stM$ is a closed embedded submanifold of $\Rnp_*$ with dimension $np-{p(p+1)}/{2}$ when $M\in\mathbb{R}^{p\times p}$ is a symmetric positive-definite matrix.
	\end{prop}
	\begin{pf}
		In the view of $\revise{\Phi_M}$, since $\stiefel$ is a closed set, we have that $\stM$ is also closed.
		{Moreover, if $\stiefel$ is locally a $\phi$-coordinate slice, then $\stM$ is locally a $\phi\circ\revise{\Phi_M^{-1}}$-coordinate slice; hence the submanifold property} 
		is preserved and by \citet[Prop. 3.3.2]{AMS08}, the set $\stM$ is an embedded submanifold with \revise{dimension} $\dimension(\stM)=\dimension(\stiefel)=np-{p(p+1)}/{2}$. 
		\qed
	\end{pf}

	\subsection{Riemannian geometry of $\stM$}
	
	We first characterize the tangent space of $\stM$.
	\begin{prop}
		The tangent space of $\stM$ at $Y\in \stM$ can be parameterized in the following ways
		\begin{subequations}
			\begin{align}
				&\TY = \{\xi\in\Rnp: \xi^{\top}Y + Y^{\top}\xi =0\} \label{eq:tangent-0}
				\\& = \{Y(Y^{\top}Y)^{-1}\varOmega+Y_{\perp}K: \varOmega\in\skewset^p, K\in\mathbb{R}^{(n-p)\times p}\}    \label{eq:tangent-1}\\
				&= \{WY: W\in\skewset^n\} \label{eq:tangent-2}\\
				&= \{\revise{\Phi_M}(\zeta): \zeta \in \mathrm{T}_{\revise{\Phi_M^{-1}}(Y)}\stiefel\}. \label{eq:tangent-3}
			\end{align}
		\end{subequations}
	\end{prop}
	\begin{pf}
		The first form can be obtained as in~\citet[\S3.3.2]{AMS08}, where $I_p$ has to be replaced by $M$ and $X\widehat{Z}$ by $YM^{-1}\widehat{Z}$.
 
		The second formulation~\eqref{eq:tangent-1} can be verified by plugging into~\eqref{eq:tangent-0} and from the fact that its dimension (${p(p-1)}/{2}+(n-p)p$) agrees with $\TY$.
		 
		The third formulation~\eqref{eq:tangent-2} can be also verified by~\eqref{eq:tangent-0} and by its dimension. 
		{Specifically, we have
			\begin{align*}
				\{WY: W\in\skewset^n\}\subseteq\TY. 
			\end{align*}
		Since $Y\in\Rnp_*$, there exists an orthogonal matrix $P\in\Rnn$ such that $P\zz Y=[I_p~0]\zz$. Let $B:= P\zz W P=\left[
		\begin{smallmatrix}
			B_{11} & B_{12}\\
			B_{21} &B_{22}\\
		\end{smallmatrix}
		\right]$. It turns out that 
		\begin{align*}
			&\dimension\{WY: W\in\skewset^n\} \\
			=\,&\dimension\{P\zz WP P\zz Y: W\in\skewset^n\} \\
			  =\,&\dimension\{B \left[
			\begin{smallmatrix}
				I_{p}\\
				0
			\end{smallmatrix}
			\right] : W\in\skewset^n\}\}  \\
			 =\,&\dimension\{\left[
			\begin{smallmatrix}
				B_{11}\\
				{B_{21}}\\
			\end{smallmatrix}
			\right]: B_{11}\in\skewset^{p}, B_{21}\in\mathbb{R}^{(n-p)\times p}\} \\
			=\,& \frac12 p(p-1) +(n-p)p = np-\frac12 p(p+1),
		\end{align*}
		which agrees with the dimension of $\TY$. The first equality comes from the fact that multiplying a subspace by an invertible matrix does not change its dimension.
		}
		
		
		The final formulation~\eqref{eq:tangent-3} follows from~\eqref{eq:tangent-0}.
		\hfill$\qed$
	\end{pf}
	Note that~\eqref{eq:tangent-2} is an over-parameterization of $\TY$ since $\dimension(\skewset^n) \,{=}\, {n(n-1)}/{2}$ and the dimension of the tangent space is only ${p(p-1)}/{2}+(n-p)p$.
	
	Given $Y\in\stM$, let $X=\revise{\Phi_M^{-1}}(Y)=YM^{-\frac12}\in\stiefel$. By \revise{making a pullback for} the canonical metric~\eqref{eq:metric-canonical}, we construct the following metric on $\Rnp_*$ for all $\xi,\zeta\in\Rnp$,
	\begin{align}
		g_Y(\xi,\zeta) &:= g^{\mathrm{c}}_{\revise{\Phi_M^{-1}}(Y)}\dkh{\revise{\Phi_M^{-1}}(\xi),\revise{\Phi_M^{-1}}(\zeta)} \nonumber\\
		&~=\langle \xi , (I_n-\frac12 YM^{-1}Y\zz) \zeta M^{-1}\rangle \nonumber\\
		&~=\langle \xi , (I_n-\frac12 Y(Y\zz Y)^{-1}Y\zz) \zeta (Y\zz Y)^{-1}\rangle.  \label{eq:metric}
	\end{align}
	Since $\revise{\Phi_M}$ is a diffeomorphism of $\Rnp_*$ onto itself and $g^{\mathrm{c}}$ is a well-defined Riemannian metric, it follows that $g$ is also a well-defined Riemannian metric on $\Rnp_*$. Hence, $(\stM,g)$ is a Riemannian submanifold of $(\Rnp_*,g)$. In particular, the metric $g$ reduces to the canonical metric when $M\,{=}\,I_p$, which implies that we generalize the canonical metric to all manifolds $\stM$. Furthermore, by construction, $\revise{\Phi_M}$ is actually an isometry between the manifolds $(\stiefel,g^{\mathrm{c}})$ and $(\stM,g)$.
	
	
	The normal space at $Y\in\stM$ with respect to $g$ is defined as the set of matrices $N\in\Rnp$ such that 
	$$g_Y(N,\xi)=0\quad\mbox{~for all~} \xi\in\TY.$$
	The following proposition gives the explicit form of any normal vector.
	\begin{prop}\label{prop:stMnormal}
		The normal space of $(\stM,g)$ at $Y\,{\in}\, \stM$ admits the following expression,
		\begin{equation}\label{eq:normal}
			\NY := \{Y(Y\zz Y)^{-1}S: S\in\symset^p\}.
		\end{equation}
	\end{prop}
	\begin{pf}
		For any $S\in\symset^p$ and $W\in\skewset^n$, it holds that 
		\begin{align*}
			&g_Y(WY,Y(Y\zz Y)^{-1}S) \\
			&=\langle WY, (I_n-\frac12 Y(Y\zz Y)^{-1}Y\zz) Y \left(Y\zz Y\right)^{-1}S (Y\zz Y)^{-1}\rangle\\
			&=  \langle  Y\zz WY, \frac12 (Y\zz Y)^{-1} S (Y\zz Y)^{-1}\rangle =0, 
		\end{align*}
		where the first equality is the definition of the metric $g$, the second equality can be verified by expanding the sum with the identity and one of the inverses canceling the $Y\zz Y$ term, and the last equality follows from $Y\zz W Y$ being skew-symmetric and therefore being orthogonal to the symmetric matrix on the right side of the inner product. According to~\eqref{eq:tangent-2}, and counting the dimension of $\TY$ and $\{Y(Y\zz Y)^{-1}S: S\in\symset^p\}$  (i.e., ${p(p+1)}/{2}$), it yields the result.
		\hfill$\qed$
	\end{pf}
	
	\subsection{Riemannian gradient of $f$ on $(\stX,g)$}
	Let $X\in\Rnp_*$. By definition of $\stM$ and the fact that $X$ is of full rank, we have that $X\in\stX$. 
	Recall \revise{that} the Riemannian gradient of $f$ with respect to the metric~$g$, denoted as $\mathrm{grad}f(X)$, is defined to be the element of $\TX$ such that
	\begin{equation}\label{eq:Riemannian_gradient}
		g_X\left(\mathrm{grad}f(X),\xi\right)=\mathrm{D}f(X)[\xi]=\langle \nabla{f}(X), \xi\rangle
	\end{equation}
	holds for all $\xi\in\TX$. The definition allows us to derive the following proposition giving a Riemannian interpretation to the relative gradient $\psi(X)X$ in the landing field~\eqref{eq:landing-field}.
	\begin{prop} \label{prop:gradient}
		The Riemannian gradient of a function $f$ on $(\stX,g)$ has the form
		\begin{equation*}\label{eq:R-gradient}
			\mathrm{grad}f(X) = \psi(X)X, 
		\end{equation*}
		where $\psi(X)\,{=}\,2\skewsym(\nabla f(X)X^{\top})$.
	\end{prop}
	\begin{pf}
		Let $W\in\skewset^n$ be any skew-symmetric matrix. In view of~\eqref{eq:tangent-2}, $WX$ is an arbitrary element of the tangent space $\TX$.  We have
		\begin{align*}
			&g_X(\psi(X)X,WX) \\
			&= \langle \psi(X)X, (I_n-\frac12 X(X\zz X)^{-1}X\zz) WX (X\zz X)^{-1}\rangle\\
			&= \langle (I_n-\frac12 X(X\zz X)^{-1}X\zz) \psi(X)X (X\zz X)^{-1} X\zz, W\rangle \\
			&= \langle \nabla{f}(X)X\zz - \sym(X(X\zz X)^{-1}X\zz \nabla f(X) X\zz), W\rangle \\
			&= \langle \nabla{f}(X)X\zz,W\rangle \\
			&= \langle \nabla{f}(X),WX\rangle, 
		\end{align*}
		where in the second equality we rearranged the terms in the inner product, the third equality can be verified by expanding the sum with the identity and writing out the definition of $\psi(X)$, the fourth equality comes from the fact that symmetric matrices are orthogonal to skew-symmetric matrices, and in the last fifth equality we move $X\zz$ to the right side of the \revise{Frobenius} inner product. According to the definition~\eqref{eq:Riemannian_gradient}, it yields the result.
		\hfill$\qed$
	\end{pf}
	
\begin{figure}[htbp]
		\centering
		\small
		\begin{tikzpicture}[scale=.4]
			\filldraw[color=cyan!10] (1.5,6) -- (-1,4) -- (9,3) -- (11,5.3) -- (1.5,6);
			
			\coordinate [label=left:$\Rnp$] (R) at (13,-4);
			\coordinate [label=-45:$X$] (X) at (5,4.8);
			\coordinate [label=left:$\stX$] (M) at (13,1);
			\coordinate [label=left:$\TX$] (T) at (13,6);
			\coordinate [label=right:$\NX$] (N) at ($(X)+(0,3)$);
			\coordinate (Y) at (2,7);  
			\coordinate [label=180:{$-\psi(X)X=-\mathrm{grad}f(X)$}] (F) at ($(Y)+(0,-2)$);
			\coordinate [label=0:{$-\lambda\nabla\mathcal{N}(X)$}] (F1) at ($(X)+(0,-2.77)$);
			\coordinate [label=-90:{$-\Lambda(X)$}] (Y1) at ($(F)-(X)+(F1)$);
			
			\draw (X) -- (N);
			\draw[dashed] (X) -- ($(X)+(0,-2.8)$);
			\draw ($(X)+(0,-2.95)$) -- ($(X)+(0,-3.5)$);
			\draw[-{stealth[red]},red,thick] (X) --  (F);
			\draw[-{stealth[blue]},blue,thick,dashed] (X) --  (F1);
			\draw[dashed] (Y1) -- (F1);
			\draw[dashed] (F) -- (Y1);
			\draw[cyan,ultra thick] (Y1) -- ($0.6*(Y1)+0.4*(X)$);
			\draw[-{stealth[cyan]},cyan,ultra thick,dashed] (X) -- (Y1);
			
			\node [fill=black,inner sep=.8pt,circle] at (X) {};
			
			\draw (0,2)	.. controls (0.5,5) and (3,5) .. (6,0.5);
			\draw (0,2)	.. controls (1,2.5) and (2,3) .. (3.92,3);
			\draw (6,0.5)	.. controls (7.5,1.5) and (9.5,2) .. (12,2);
			\draw (0.5,3.4)	.. controls (1.7,5.5) and (11,8) .. (12,2);
			
			\draw [densely dotted] (1.5,6) -- (-1,4) -- (9,3) -- (11,5.3) -- (1.5,6);			
			
			\coordinate (P) at (-6,-5); 
			
			\draw ($(0,2)+(P)$)	.. controls ($(0.5,5)+(P)$) and ($(3,5)+(P)$) .. ($(6,0.5)+(P)$);
			\draw ($(0,2)+(P)$)	.. controls ($(1,2.5)+(P)$) and ($(2,3)+(P)$) .. ($(3.92,3)+(P)$);
			\draw ($(6,0.5)+(P)$)	.. controls ($(7.5,1.5)+(P)$) and ($(9.5,2)+(P)$) .. ($(12,2)+(P)$);
			\draw ($(0.5,3.4)+(P)$)	.. controls ($(1.7,5.5)+(P)$) and ($(11,8)+(P)$) .. ($(12,2)+(P)$);
			
			\draw[dashed,rotate=20] ($(8.3,0.8)+(P)$) ellipse (5mm and 2.5mm);
			\draw[dashed,rotate=20] ($(8.3,0.8)+(P)$) ellipse (15mm and 7.5mm);
			\draw[dashed,rotate=20] ($(8.3,0.8)+(P)$ )ellipse (25mm and 12.5mm);
			\draw[dashed] ($(4,3)+(P)$)	.. controls ($(4.5,4.5)+(P)$) and ($(7,6)+(P)$) .. ($(9,5.4)+(P)$);
			\draw[dashed] ($(5.3,1.5)+(P)$)	.. controls ($(7,1.5)+(P)$) and ($(9,1)+(P)$) .. ($(11,4.2)+(P)$);
			
			\node at ($(10.9,2.5)+(P)$) {$f$};
			
			\coordinate [label=left:${\stiefel=\revise{\mathrm{St}_{I_p}(p,n)}}$] (M1) at ($(15,1)+(P)$);
			
		\end{tikzpicture}
		\caption{Geometric interpretation of the orthogonal components of the landing field. \label{fig:interpretation}}
	\end{figure}
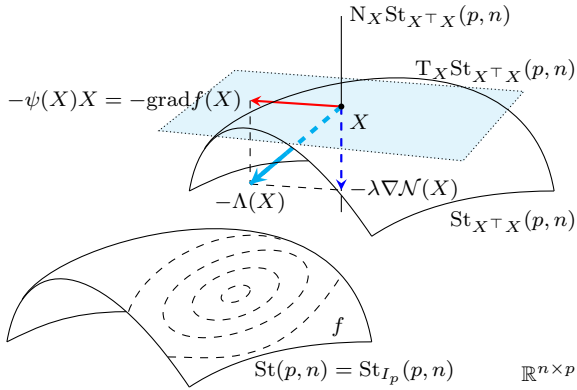

	The above Proposition~\ref{prop:gradient} and the expression of the normal space in Proposition~\ref{prop:stMnormal} give a clear interpretation of both components of the landing field $\Lambda(X) = \psi(X)X+\lambda \nabla\mathcal{N}(X)$. Specifically, $\psi(X)X$ is the Riemannian gradient of $f$ on the submanifold $(\stX,g)$, and 
	$$\nabla\mathcal{N}(X)=X(X\zz X)^{-1}((X\zz X)^2-X\zz X)$$
	belongs to the normal space $\NX$. \revise{Consequently}, $\Lambda(X)$ is the linear combination of the two orthogonal fields in the tangent and the normal space of $(\stX,g)$; see Fig.~\ref{fig:interpretation} for a geometric illustration. The orthogonal property will have important consequences in the next section where we analyze convergence of the landing flow.

	\section{Convergence of the landing flow} \label{sec:convergence}
	In this section, we establish a convergence \revise{analysis} for the solutions of the landing system~\eqref{eq:landing-flow}, denoted as $\varphi_t(X_0)$ for a starting point $X_0 \in\Rnp_*$ and for all $t\ge0$. The proof consists of two parts, firstly by the convergence of $X(t)\zz X(t)$ to $I_p$ and secondly by the convergence of $X(t)$ to the set of critical points of $f$ relative to $\stiefel$.
	
	{{\it Standing assumption:} $\nabla f$ is locally Lipschitz continuous.}
	
	We show that the solutions of the landing system exist and are unique, thus making the landing flow well defined.
	\begin{prop}[Existence and uniqueness]\label{prop:existence}
		For the landing system~\eqref{eq:landing-flow} starting at $X_0\in\Rnp_*$, there exists a unique solution $t\mapsto \flow$ defined for all $t \geq 0$ such that $\varphi_0(X_0)=X_0$. Moreover, we have that  $\mathcal{N}(\flow)$ is nonincreasing.
	\end{prop}
	\begin{pf}
		Differentiating $\mathcal{N}(X(t))$ with respect to $t$ gives
		\begin{align*}
			\dt\mathcal{N}(X(t))&=\langle{\dot{X}(t),\nabla \mathcal{N}(X(t))}\rangle\\
			&=-\langle{\psi(X(t))X(t)+\lambda \nabla \mathcal{N}(X(t)), \nabla \mathcal{N}(X(t))}\rangle \\
			&= -\lambda \norm{\nabla \mathcal{N}(X(t))}\fs \leq0,
		\end{align*}
		where the second equality comes from the definition of the landing field and the last equality is the consequence of $\psi(X)X$ being orthogonal to $\nabla \mathcal{N}(X(t))$. Hence, $\mathcal{N}(X(t))$ is nonincreasing and each solution of the landing system remains in a compact set. By $\nabla f$ being locally Lipschitz, we have that $\Lambda(X)$ is also locally Lipschitz. By the Picard–Lindelöf theorem, the landing system has a unique solution.
		\hfill$\qed$
	\end{pf}
	It is worth noting that Proposition~\ref{prop:existence} holds for any $\lambda>0$. This is due to the orthogonality of the two components of the landing field. By contrast with the landing field, the components of the PLAM field $\tilde{\Lambda}$ defined in~\eqref{eq:PLAM-field} do not satisfy the orthogonal property, and as a consequence, the existence of its flow requires a lower threshold on $\lambda>\lambda_0>0$ \citep{gao2018parallelizable}.
	
	{Since $\mathcal{N}(\flow)$ is nonincreasing and the set of minimizers of $\mathcal{N}$ is the Stiefel manifold, it follows that the Stiefel manifold is an invariant of the landing flow. Recall also that, on the Stiefel manifold, the landing flow reduces to the Riemannian gradient flow with respect to the canonical metric.}
	
	\subsection{Convergence of $X\zz X$ to $I_p$}
	The following result shows that the landing flow $\varphi_t(X_0)$ converges to $\stiefel$ as $t\to\infty$ for any $X_0\in\Rnp_*$.
	\begin{prop}[Convergence to the Stiefel manifold]\label{prop:convergence-stiefel}
		For \linebreak all $X_0\in\Rnp_*$, 
		{we have that $\varphi_t(X_0) \in \mathbb{R}_*^{n\times p}$ for all $t>0$ and}
		\begin{equation*}
			\lim_{t\to\infty} \mathcal{N}(\flow) = 0.
		\end{equation*}
	\end{prop}
	\begin{pf}
		Let $\chi(t) := X(t)\zz X(t)$ with $X(t)$ following the dynamics of the landing system~\eqref{eq:landing-flow}. Differentiating $\chi(t)$ with respect to $t$ yields
		\begin{align*}
			\dot \chi(t) &= \dot{X}(t)\zz X(t) + X(t)\zz \dot{X}(t) \\
			&= -2\lambda\, \chi(t)\left(\chi(t)-I_p\right)
		\end{align*}
		By the right hand side being a matrix polynomial function of a symmetric matrix $\chi(t)$, we have that $\chi(t)$ has constant eigenvectors for all $t\geq0$ and its eigenvalues $\{\chi_i\}_{i=1}^p$ follow $\dot \chi_i(t) = -2\lambda\, \chi_i(t) (\chi_i(t)-1)$. The solution of the ODE for the eigenvalues can be computed explicitly as 
		\begin{equation*}
			\chi_i(t)=\frac{\chi_i(0)e^{2\lambda t}}{\chi_i(0)(e^{2\lambda t}-1)+1}.
		\end{equation*}
		Since $\lambda>0$ and $\chi_i(0)>0$ because $X_0\in\Rnp_*$ is of full rank, we have that $\lim_{t\to\infty}\chi_i(t)=1$, i.e., all eigenvalues of $\chi(t)$ converge to~$1$. Hence $\chi(t)$ converges to $I_p$, and thus $\mathcal{N}(\flow)$ converges to~$0$.
		\hfill$\qed$
	\end{pf}
	
	\subsection{Convergence of the landing flow}
	Let $\Zstar$ denote the set of critical points of $f$ relative to $\stiefel$. Since $\psi(X)X$ is the Riemannian gradient on $\stiefel$ with the canonical metric when $X\in\stiefel$, according to~\citet[\S4.1]{AMS08}, it follows that
	\begin{equation*}\label{eq:critical-set}
		\Zstar=\left\{X^*\in\stiefel: \psi(X^*)X^*=0\right\}.
	\end{equation*}
	Moreover, by the orthogonality of the two terms in the landing field $\Lambda$, we have that
	\begin{equation*}
		X^*\in\Zstar\quad \mbox{if and only if}\quad \Lambda(X^*)=0.    
	\end{equation*}
	
	Recall that the $\omega$-limit set of a trajectory $\varphi_t$ is the set of points $\varphi^*$ for which there exists a sequence $\{t_n\}$ with $\lim_{n\to \infty} t_n=\infty$ such that $\lim_{n\to \infty} \varphi_{t_n}=\varphi^*$; see e.g.,~\cite{Kha96}.
	
	Next, we show that the landing flow converges to the set of critical points of $f$ on $\stiefel$.
	
	\begin{thm}[Convergence to the set of critical points]\label{thm:omega-limit}
		For all $X_0\in\Rnp_*$, the $\omega$-limit points of $\flow$ belong to~$\Zstar$; in other words, the solution of the landing system~\eqref{eq:landing-flow} converges to the set of critical points of $f$ relative to the Stiefel manifold $\mathrm{St}(n,p)$.
	\end{thm}
	\begin{pf}
		Let $b\geq a\geq0$ and 
		\[
		\mathcal{N}^{-1}([a,b]) := \{X\in\Rnp_*: a \leq \mathcal{N}(X) \leq b\}.
		\]
		Let $\epsilon>0$ and 
		\[
		B_\epsilon(\Zstar)
		:= \bigcup_{{X^{*}}\in\Zstar} B_\epsilon(X^*)
		\]
		where $B_\epsilon(X^*) := \left\{X\in\Rnp_*: \left\|X-X^*\right\|\ff < \epsilon\right\}$.
		
		{We claim (to be proven in the next paragraph) that for any}
		$\epsilon>0$ there exists $\delta>0$ such that
		\begin{equation*}
			\max_{ X\in \mathcal{N}^{-1}([0,\delta]) \setminus B_\epsilon(\Zstar) } \mathrm{D} f(X)[-\Lambda(X)] < 0.
		\end{equation*} 
		According to Proposition~\ref{prop:convergence-stiefel}, the solution $\flow$ eventually stays in $\mathcal{N}^{-1}([0,\delta])$. As a result, $\flow$ converges to $B_\epsilon(\Zstar)$; otherwise, due to the above claim, $\lim_{t\to\infty}f(\flow)=-\infty$, which is impossible since $f$ is continuous and $\mathcal{N}^{-1}([0,\delta])$ is compact. Since the convergence of $\flow$ to $B_\epsilon(\Zstar)$ holds for all $\epsilon>0$, it follows that $\flow$ converges to $\Zstar$.
		
		We now show by contradiction that the claim in the first part of the proof is true. Suppose the statement is not true, that is, there exists $\delta_k>0$ monotononically decreasing with $\lim_{k\to\infty}\delta_k=0$ and $X_k \in \mathcal{N}^{-1}([0,\delta_k]) \setminus B_\epsilon(\Zstar)$ such that 
		\begin{equation*}\label{eq:omega-limit1}
			\mathrm{D} f(X_k)[-\Lambda(X_k)] \geq 0.
		\end{equation*}
		Since $\mathcal{N}^{-1}([0,\delta_0]) \setminus B_\epsilon(\Zstar)$ is compact, the sequence $\{X_k\}$ has a convergent subsequence. Let $\tilde{X}$ be its limit. On the one hand, by continuity of $X\mapsto \mathrm{D} f(X)[-\Lambda(X)]$, we have $\mathrm{D} f(\tilde{X})[-\Lambda(\tilde{X})] \geq 0$. On the other hand, $\tilde{X} \in \stiefel$ since $\mathcal{N}(X_k)\leq \delta_k$ and $\lim_{k\to\infty}\delta_k=0$, and moreover $\tilde{X}\notin\Zstar$ since it is at least a distance $\epsilon$ away from the critical points. \revise{This further implies that $\Lambda(\tilde{X})=\psi(\tilde{X})\tilde{X}\neq0$ is the relative gradient of $f$ on $\stiefel$.}
		Hence we have $\mathrm{D} f(\tilde{X})[-\Lambda(\tilde{X})] < 0$, a contradiction.
		\hfill$\qed$
	\end{pf}
	
	\subsection{Stability of the equilibria}
	We investigate the stability of the equilibria of the landing system~\eqref{eq:landing-flow} for minimizing $f$ relative to $\stiefel$. Note that \cite{absil2004continuous} considered the continuous-time flows on quotient spaces by using a similar idea. 
	
	By definition (see e.g.,~\cite{absil2006on}), $X^*$ is an equilibrium point of the system~\eqref{eq:landing-flow} if $\Lambda(X^*)=0$, i.e., $X^*$ is a critical point of $f$ relative to $\stiefel$. Furthermore, $X^*$ is stable if, for any $\epsilon>0$, there exists $\delta>0$ such that, $\norm{X(0)-X^*}<\delta$ implies $\norm{X(t)-X^*}<\epsilon$ for all $t\geq0$. Moreover, it is asymptotically stable if it is stable and
	there exists $\delta>0$ such that, $\norm{X(0)-X^*}<\delta$ implies $\lim_{t\to\infty} X(t)=X^*$.
	
	Next, we give a sufficient condition for asymptotic stability. The proof is based on the theory of semidefinite Lyapunov functions; see~\cite{iggidr1996semidefinite}. 
	\begin{thm}[Asymptotic stability]\label{thm:stability}
		If $X^*$ is a local minimum and isolated critical point of $f$ relative to $\stiefel$, then $X^*$ is an asymptotically stable point of~the landing system~\eqref{eq:landing-flow}. 
	\end{thm}
	\begin{pf}
		\revise{Let $\dot{\mathcal{N}}(X)$ denote $\dt \mathcal{N}(\varphi_t(X))|_{t=0}$.}
		
		(i) Since $X^*$ is a critical point, it is also an equilibrium point of~\eqref{eq:landing-flow}, and $X^*\in\stiefel$, i.e., $\mathcal{N}(X^*)=0$.  In addition, we have \revise{ $\mathcal{N}(X)\geq 0$ for all $X$}. 
		
		(ii) According to the proof of Proposition~\ref{prop:existence}, it holds that \revise{$\dot{\mathcal{N}}(X)\leq 0$ for all $X \in \Rnp_*$}.
		
		\revise{(iii) We have $$\{X\in\Rnp_*: \dot{\mathcal{N}}(X)=0\}=\stiefel,$$
		which is an invariant of the landing flow.
		
		(iv) Since the landing flow is a gradient descent flow for $f$ relative to $\stiefel$, and since $X^*$ is a local minimum and isolated critical point of $f$ relative to $\stiefel$, it follows from~\citet[\S4]{absil2006on} that $X^*$ is asymptotically stable relative to $\stiefel$.
		
		
		The above points combined with~\citet[Corollary~1]{iggidr1996semidefinite} yield the result that $X^*$ is an asymptotically stable equilibrium point.
		}
		\hfill$\qed$
	\end{pf}
	\begin{cor}
		For all $X_0\in\Rnp_*$, if $X^*$ is a local minimum and isolated critical point of $f$ relative to $\stiefel$, and if $X^*$ is an $\omega$-limit point of $\flow$, then $\lim_{t\to\infty} \flow = X^*$.
	\end{cor}
	\begin{pf}
		Since $X^*$ is an $\omega$-limit point, $\varphi_t(X_0)$ eventually enters any neighborhood of $X^*$. Since moreover, in view of Theorem~\ref{thm:stability}, $X^*$ is (asymptotically) stable, it follows that $\varphi_t(X_0)$ eventually stays in any neighborhood of $X^*$.
		\hfill$\qed$
	\end{pf}
	
	
	\section{Conclusion}\label{sec:conclusion}
	We have proposed an extension of the landing flow of~\cite{ablin2021fast} to rectangular matrices, obtained a Riemannian gradient interpretation to the $\psi(X)X$ term of the landing field~\eqref{eq:landing-field}, and proven that the solutions of the landing system globally converge to the set of equilibria of the objective function relative to the Stiefel manifold. In future work, we will address the question of finding a discrete-time counterpart of the landing flow that preserves its favorable convergence properties. 
	
	%

\begin{thebibliography}{22}
	\providecommand{\natexlab}[1]{#1}
	\providecommand{\url}[1]{\texttt{#1}}
	\providecommand{\urlprefix}{URL }
	\expandafter\ifx\csname urlstyle\endcsname\relax
	\providecommand{\doi}[1]{doi:\discretionary{}{}{}#1}\else
	\providecommand{\doi}{doi:\discretionary{}{}{}\begingroup
		\urlstyle{rm}\Url}\fi
	
	\bibitem[{Ablin and Peyr\'e(2022)}]{ablin2021fast}
	Ablin, P. and Peyr\'e, G. (2022).
	\newblock Fast and accurate optimization on the orthogonal manifold without
	retraction.
	\emph{Proceedings of The 25th International Conference on Artificial
		Intelligence and Statistics}, volume 151 of \emph{Proceedings of Machine
		Learning Research}, 5636--5657. PMLR.
	\newblock \urlprefix\url{https://proceedings.mlr.press/v151/ablin22a.html}.
	
	\bibitem[{Absil(2004)}]{absil2004continuous}
	Absil, P.-A. (2004).
	\newblock Continuous-time flows on quotient spaces for principal component
	analysis.
	\newblock In \emph{Proceedings of the 16th International Symposium on
		Mathematical Theory of Networks and Systems (MTNS2004)}.
	
	\bibitem[{Absil(2006)}]{Absil2006continuous}
	Absil, P.-A. (2006).
	\newblock Continuous-time systems that solve computational problems.
	\newblock \emph{International Journal of Unconventional Computing}, 2(4),
	291--304.
	
	\bibitem[{Absil and Kurdyka(2006)}]{absil2006on}
	Absil, P.-A. and Kurdyka, K. (2006).
	\newblock On the stable equilibrium points of gradient systems.
	\newblock \emph{Systems $\&$ Control Letter}, 55(7), 573--577.
	\newblock \doi{10.1016/j.sysconle.2006.01.002}.
	
	\bibitem[{Absil et~al.(2008)Absil, Mahony, and Sepulchre}]{AMS08}
	Absil, P.-A., Mahony, R., and Sepulchre, R. (2008).
	\newblock \emph{Optimization Algorithms on Matrix Manifolds}.
	\newblock Princeton University Press.
	
	\bibitem[{Boumal(2022)}]{boumal2022intromanifolds}
	Boumal, N. (2022).
	\newblock \emph{An introduction to optimization on smooth manifolds}.
	\newblock To appear with Cambridge University Press.
	\newblock \urlprefix\url{http://www.nicolasboumal.net/book}.
	
	\bibitem[{Bournez and Pouly(2021)}]{bournez2021survey}
	Bournez, O. and Pouly, A. (2021).
	\newblock A survey on analog models of computation.
	\newblock In \emph{Handbook of Computability and Complexity in Analysis},
	173--226. Springer.
	\newblock \doi{10.1007/978-3-030-59234-9_6}.
	
	\bibitem[{Brockett(1991)}]{Brockett1991dynamical}
	Brockett, R.W. (1991).
	\newblock Dynamical systems that sort lists, diagonalize matrices, and solve
	linear programming problems.
	\newblock \emph{Linear Algebra and its Applications}, 146, 79--91.
	\newblock \doi{10.1016/0024-3795(91)90021-N}
	
	\bibitem[{Chen et~al.(2020)Chen, Ma, {Man-Cho So}, and
		Zhang}]{Chen2020Proximal}
	Chen, S., Ma, S., {Man-Cho So}, A., and Zhang, T. (2020).
	\newblock Proximal gradient method for nonsmooth optimization over the
	Stiefel manifold.
	\newblock \emph{SIAM Journal on Optimization}, 30(1), 210--239.
	\newblock \doi{10.1137/18M122457X}.
	
	\bibitem[{Chu(1994)}]{Chu1994list}
	Chu, M.T. (1994).
	\newblock A list of matrix flows with applications.
	\newblock In \emph{Hamiltonian and gradient flows, algorithms and control},
	volume~3 of \emph{Fields Institute Communications}, 87--97. 
	
	\bibitem[{Edelman et~al.(1998)Edelman, Arias, and Smith}]{edelman1998geometry}
	Edelman, A., Arias, T.A., and Smith, S.T. (1998).
	\newblock The geometry of algorithms with orthogonality constraints.
	\newblock \emph{SIAM Journal on Matrix Analysis and Applications}, 20(2), 303--353.
	\newblock \doi{10.1137/S0895479895290954}.
	
	\bibitem[{Eld{\'e}n and Park(1999)}]{Elden1999Procrustes}
	Eld{\'e}n, L. and Park, H. (1999).
	\newblock A {{Procrustes}} problem on the {{Stiefel}} manifold.
	\newblock \emph{Numerische Mathematik}, 82(4), 599--619.
	\newblock \doi{10.1007/s002110050432}.
	
	\bibitem[{Gao et~al.(2022)Gao, Hu, Kuang, and Liu}]{gao2020orthogonalization}
	Gao, B., Hu, G., Kuang, Y., and Liu, X. (2022).
	\newblock An orthogonalization-free parallelizable framework for all-electron
	calculations in density functional theory.
	\newblock \emph{SIAM Journal on Scientific Computing}, 44(3), B723--B745.
	\newblock \doi{10.1137/20M1355884}.
	
	\bibitem[{Gao et~al.(2019)Gao, Liu, and Yuan}]{gao2018parallelizable}
	Gao, B., Liu, X., and Yuan, Y.-X. (2019).
	\newblock Parallelizable algorithms for optimization problems with
	orthogonality constraints.
	\newblock \emph{SIAM Journal on Scientific Computing}, 41(3), A1949--A1983.
	\newblock \doi{10.1137/18M1221679}.
	
	\bibitem[{Golub and Van~Loan(2013)}]{golub2013matrix}
	Golub, G.H. and Van~Loan, C.F. (2013).
	\newblock \emph{Matrix Computations}.
	\newblock Johns Hopkins University Press, 4th edition.
	
	\bibitem[{Grubi{\v s}i{\'c} and Pietersz(2007)}]{Grubisic2007Efficient}
	Grubi{\v s}i{\'c}, I. and Pietersz, R. (2007).
	\newblock Efficient rank reduction of correlation matrices.
	\newblock \emph{Linear Algebra and its Applications}, 422(2-3), 629--653.
	\newblock \doi{10.1016/j.laa.2006.11.024}.
	
	\bibitem[{Hu et~al.(2020)Hu, Liu, Wen, and Yuan}]{Hu2020Brief}
	Hu, J., Liu, X., Wen, Z.W., and Yuan, Y.-X. (2020).
	\newblock A brief introduction to manifold optimization.
	\newblock \emph{Journal of the Operations Research Society of China}, 8(2),
	199--248.
	\newblock \doi{10.1007/s40305-020-00295-9}.
	
	\bibitem[{Iggidr et~al.(1996)Iggidr, Kalitine, and
		Outbib}]{iggidr1996semidefinite}
	Iggidr, A., Kalitine, B., and Outbib, R. (1996).
	\newblock Semidefinite Lyapunov functions stability and stabilization.
	\newblock \emph{Mathematics of Control, Signals and Systems}, 9(2), 95--106.
	\newblock \doi{10.1007/BF01211748}.
	
	\bibitem[{Joho and Mathis(2002)}]{Joho2002Joint}
	Joho, M. and Mathis, H. (2002).
	\newblock Joint diagonalization of correlation matrices by using gradient
	methods with application to blind signal separation.
	\newblock In \emph{Sensor {{Array}} and {{Multichannel Signal Processing
				Workshop Proceedings}}, 2002}, 273--277.
	\newblock \doi{10.1109/SAM.2002.1191043}.
	
	\bibitem[{Khalil(1996)}]{Kha96}
	Khalil, H.K. (1996).
	\newblock \emph{Nonlinear systems}.
	\newblock Prentice Hall, second edition.
	
	\bibitem[{Mahony et~al.(1996)Mahony, Helmke, and Moore}]{MHM96}
	Mahony, R.E., Helmke, U., and Moore, J.B. (1996).
	\newblock Gradient algorithms for principal component analysis.
	\newblock \emph{The Journal of the Australian Mathematical Society. Series B. Applied Mathematics}, 37(4), 430--450.
	\newblock \doi{10.1017/S033427000001078X}
	
	\bibitem[{Wen and Yin(2013)}]{Wen2013feasible}
	Wen, Z. and Yin, W. (2013).
	\newblock A feasible method for optimization with orthogonality constraints.
	\newblock \emph{Mathematical Programming}, 142(1-2), 397--434.
	\newblock \doi{10.1007/s10107-012-0584-1}.
	
\end{thebibliography}

\end{document}